\font\script=eusm10.
\font\sets=msbm10.
\font\stampatello=cmcsc10.
\def\0{{\bf 0}}
\def\1{{\bf 1}}
\def\defineq{\buildrel{def}\over{=}}
\def\definiz{\buildrel{def}\over{\Longleftrightarrow}}

\def\C{\hbox{\sets C}}
\def\N{\hbox{\sets N}}

\def\P{\hbox{\sets P}}

\def\cloud0{\hbox{$<\0>$}}
\def\corsivoF{\hbox{\script F}}

\def\corsivoN{\hbox{\script N}}

\def\corsivoR{\hbox{\script R}}

\def\square{\hbox{\vrule\vbox{\hrule\phantom{s}\hrule}\vrule}}

\def\qed{\hfill $\square$}
\def\inzero{\in\,<\0>}

\par
\centerline{\bf Multiplicative Ramanujan coefficients of null-function}
\bigskip
\bigskip
\centerline{Giovanni Coppola and Luca Ghidelli}\footnote{ }{MSC$2010$: $11{\rm A}25$, $11{\rm K}65$, $11{\rm N}37$ - Keywords: Ramanujan expansion, Ramanujan coefficient, finite Euler product, absolute convergence, null-function} 
\bigskip
\bigskip
\par
\noindent
{\bf Abstract}. The null-function $\0(a)\defineq 0$ $\forall a \in \N$ has classical {\it Ramanujan expansions}: $\0(a)=\sum_{q=1}^{\infty}(1/q)c_q(a)$ (here $c_q(a)$ is the {\it Ramanujan sum}), given by Ramanujan himself (in his 1918 history-making paper), and $\0(a)=\sum_{q=1}^{\infty}(1/\varphi(q))c_q(a)$, given by Hardy few years later ($\varphi$ is Euler's totient function). Both are pointwise converging in any $a\in \N$, but not absolutely convergent. A general $G:\N \rightarrow \C$ is called a {\it Ramanujan coefficient}, abbrev. R.c., iff (if and only if) $\sum_{q=1}^{\infty}G(q)c_q(a)$ converges in all $a\in \N$; also, for general $F:\N \rightarrow \C$ we call $<F>$, the set of its R.c.s, the {\it Ramanujan cloud} of $F$. Our Main Theorem in {\tt arxiv:1910.14640}, for Ramanujan expansions and {\it finite Euler products}, applies to give a complete Classification for multiplicative Ramanujan coefficients of $\0$. We find that Ramanujan's coefficient of $\0$, $G_R(q):=1/q$, is a {\it normal} arithmetic function $G$, i.e., multiplicative with $G(p)\neq 1$ on all primes $p$; while Hardy's $G_H(q):=1/\varphi(q)$ is a {\it sporadic} $G$, namely multiplicative, $G(p)=1$ for a (necessarily finite) set of primes, but there's no prime $p$ with $G(p^K)=1$ on all integer powers $K\ge 0$ (Hardy's has \thinspace $G_H(p)=1$ \thinspace iff $p=2$). The multiplicative $G:\N \rightarrow \C$ such that there exists at least one prime $p$ with $G(p^K)=1$, on all powers $K\ge 0$, are defined to be {\it exotic}. So, this definition completes the cases for multiplicative Ramanujan coefficients of $\0$. \enspace The exotic ones are a kind of new phenomenon in the \lq \lq $\0-$cloud\rq \rq (i.e., $<\0>$): {\bf exotic Ramanujan coefficients represent } $\0$ {\bf only with a convergence hypothesis}. The not exotic, apart from the convergence of $\sum_{(q,a)=1}G(q)\mu(q)$ in any $a\in\N$, require furthermore $\sum_{q=1}^{\infty}G(q)\mu(q)=0$ for normal $G\in <\0>$, while sporadic $G\in <\0>$ need \thinspace $\sum_{(q,P(G))=1}G(q)\mu(q)=0$ (where $P(G)$ is the product of all primes $p$ making $G(p)=1$). We give many examples of multiplicative $G\in <\0>$; we prove, in passing, that the only multiplicative Ramanujan coefficients giving an absolutely converging Ramanujan expansion of $\0$ are the exotic ones; actually, these $G\in <\0>$  \thinspace generalize to the \lq \lq {\bf weakly exotic}\rq \rq \thinspace coefficients of $\0$, which {\bf are not necessarily multiplicative}. 
\bigskip
\bigskip
\bigskip

\par
\noindent{\bf Table of Contents}
\bigskip

\def\tablewidth{\hsize}
\tabskip=5em  
\halign to \tablewidth 
{
\hfil# \tabskip=2em plus10em  &
#\hfil \tabskip=7em plus 2em &
\hfill#\cr
1. & Introduction & \bf 1 \cr
2. & Statement of the classification & \bf 3 \cr
3. & Lemmata & \bf 4 \cr
4. & Proof of the Classification Theorem & \bf 6 \cr
5. & Examples of Ramanujan expansions of  $\0$ via analytic methods & \bf 8 \cr
6. & Examples of exotic and weakly-exotic Ramanujan coefficients & \bf 9 \cr
7. & Absolutely convergent Ramanujan expansions & \bf 11 \cr
}

\bigskip
\bigskip
\par
\noindent {\bf 1. Introduction}

\bigskip

\par
\noindent
Following Ramanujan [R] we denote by $c_q(a)$ the sum of $q$-th roots of unity given by (henceforth $(h,q)$ is the gcd of $h,q$) 
$$
c_q(a)\defineq \sum_{{h=1 \atop (h,q)=1}}^{q-1} e^{2\pi i a h / q} =\sum_{{h=1 \atop (h,q)=1}} ^{q-1}\cos(2\pi a h  / q)
$$
and we call it a {\it Ramanujan sum}. 
Then a {\it Ramanujan coefficient} is any arithmetic function $G:\N \rightarrow \C$ such that  
$$
\sum_{q=1}^{\infty}G(q)c_q(a)
$$
\par				
\noindent
converges pointwise for all $a\in \N$. If $F(a)= \sum_{q=1}^{\infty}G(q)c_q(a)$ we also say that the series is a {\it Ramanujan expansion} (or Fourier-Ramanujan expansion) of the arithmetic function $F:\N\rightarrow \C$, or that $G$ is a Ramanujan coefficient of $F$. 
\par

Ramanujan sums are ubiquitous expressions in number theory with remarkable properties. For example, it was first observed by Hardy that the function $c_q(a)$ is multiplicative in the variable $q$, so it is determined by its values at powers of primes $c_{p^k}(a)$. Other formul\ae \thinspace to compute the Ramanujan sums in terms of the M\"obius function $\mu$ (see [T] esp.) and Euler's totient function $\varphi$ [T] are 
$$
c_q(a) = \sum_{d\mid (q,a)} \mu\left( {q \over d} \right) d =  \mu\left( {q\over (q,a)} \right) {\varphi (q) \over \varphi(q/(q,a))}. 
$$
The first expression was found by Kluyver [K] while the second is due to von Sterneck and H\"older [H\"o]. 

\par
For a fixed $F:\N \rightarrow \C$ we denote the set of its Ramanujan coefficients by 
$$
<F>\defineq \left\{ G:\N \rightarrow \C \enspace | \enspace \forall a\in \N, \thinspace \sum_{q=1}^{\infty}G(q)c_q(a)=F(a) \right\} .
$$
\par
\noindent
This set is called the {\it Ramanujan cloud} of $F$ or, for brevity, the $F$-cloud.  
The problem of establishing the elements of  $F$-clouds (i.e., finding the arithmetic functions $G$ that are Ramanujan coefficients of $F$) has been considered by many authors. 
One general method to achieve this goal is through mean-value theorems and harmonic analysis. 
In this sense, a fundamental tool which is often exploited is an orthogonality relation for Ramanujan sums discovered by Carmichael [Ca].  For more on the classical and modern theory of Ramanujan sum and Ramanujan expansions, we refer the reader to the surveys [C1], [L1], [L2], [M] and to the monograph [ScSp]. 
We also mention that there is a theory of Ramanujan expansion for functions of several variables and of ``dual Ramanujan expansions'' of the form $\sum_{a=1}^\infty F(a)c_q(a)$ : compare [M].
\par
Another interesting problem, which is more related to the results of the present article, is that of the \lq \lq cardinality\rq \rq \thinspace of $<F>$. It was already known to Ramanujan [R] and Hardy [H] that an arithmetic function may be representable as a Ramanujan expansion in more than one way. Respectively, they found the following two nontrivial expansions of the constant zero function $\0(a)\defineq 0$, $\forall a\in \N$: 
$$
\sum_{q=1}^\infty {1 \over q} c_q(a) = 0, 
\quad\quad 
\sum_{q=1}^\infty {1\over \varphi(q)} c_q(a) = 0. 
$$
An intriguing fact is that the convergence of these series (to zero, for every $a\in\N$) turns out to be logically equivalent, in a precise sense, to the Prime Number Theorem. 
Until recently, the Ramanujan coefficients $G_R(q):=1/q$ and $G_H(q):=1/\varphi(q)$ corresponding to the two examples above of Hardy and Ramanujan were essentially the only elements of the $\0$-cloud (together of course with their linear combinations) to be found in the literature.  
Using the theory of Euler products of Ramanujan expansions the first author was able to produce several new examples of $G\inzero$ and he discovered, in particular, a whole new class of Ramanujan expansions of $\0$, which often have a much better convergence rate to zero than the examples of Hardy and Ramanujan given above. These Ramanujan expansions correspond to the class of so-called {\stampatello exotic} Ramanujan coefficients, described later in this article. 
Also, since to any Ramanujan coefficient we can add an arbitrary $G\inzero$ without changing the value of the corresponding Ramanujan expansion, it follows that every $F$-cloud has infinite cardinality, and in fact it has a natural structure of an infinite-dimensional complex affine space [C4].
\par
In this paper we shall provide a classification of the multiplicative Ramanujan coefficients $G\inzero$ under some mild convergence hypotheses which are similar to those found in work of Lucht [L2, Thm 3.1] and which come out naturally from the theory of finite Euler products of [C3]. The precise statements are given in the next section. 
With this classification at hand, in sections 5 and 6 of this article we provide many new examples of Ramanujan expansion of the constant function $\0$. 
Since the Ramanujan sums $c_q(a)$ are multiplicative themselves, the decision to classify the multiplicative portion of the $\0$-cloud is rather natural. At any rate, our results generalise easily to provide also some new examples of non-multiplicative Ramanujan coefficients $G\inzero$, such as the class of {\stampatello weakly exotic} Ramanujan coefficients described in section 6.2. 
\par				
\noindent
By the way, with the only exception of this section, our $G:\N \rightarrow \C$ will always be {\stampatello multiplicative}. 
Finally, in section 7, we shall discuss criteria of absolute convergence for Ramanujan expansions.
\par
The authors plan to give a more complete account of the non-multiplicative part of the $\0$-cloud, in a future publication. 

\bigskip

\par
\noindent {\bf 2. Statement of the classification}

\bigskip

\par
\noindent
In order to state and prove our results it is better to introduce some notation. 
We denote by $\P$ the set of prime numbers and by $\N_0\defineq \N\cup\{0\}$ the set of non-negative integers. 
For every $a\in\N$ we let $\P(a)$ be the set of prime divisors of  $a$ (recall \enspace $\P(1)\defineq \emptyset$) and we let 
(recalling the convention : empty products are $1$) 
$$
P(a)\defineq \prod_{p\in \P(a)}p 
$$
\par
\noindent
be its {\it squarefree kernel} (also known as the {\it radical} of $a$). 
By extension, for a given set \enspace $\corsivoN \subset \N$, we let \hfill $\P(\corsivoN)\defineq \{p\in \P : p|n, \hbox{\rm for\enspace some}\enspace n\in \corsivoN \}$ and then, whenever this set of primes is finite, we denote by $$P(\corsivoN)\defineq \prod_{p\in \P(\corsivoN)} p$$ their product. 
We also recall that $v_p(a)\defineq \max\{ K\in \N_0 : p^K\vert a\}$ denotes the $p-${\stampatello adic valuation} 
of a number $a\in \N$  with respect to $p\in \P$. \hfil Henceforth, we abbreviate with \lq \lq iff\rq \rq \thinspace the expression \lq \lq if and only if\rq \rq.
\par 
Now fix a multiplicative function $G:\N \rightarrow \C$. We say that $G$ is {\it multiplicatively trivial} (resp. {\it completely multiplicatively trivial}) at a prime $p\in\P$ iff $G(p)=1$ (resp. $G(p^K)=1$ for all $K\ge 0$). We can then define two kinds of ``multiplicative spectra'' of $G$ which are the key to our Classification: 
$$
\corsivoF(G)\defineq \{ p\in \P : G(p)=1 \},
\quad \quad 
\corsivoF_0(G)\defineq \{ p\in \P : G(p^K)=1, \forall K\in \N_0 \}.
$$
By multiplicativity, the function $G$ is multiplicatively (resp. completely multiplicatively) trivial at $p$ iff the equality $G(pr)=G(r)$ holds for every $(r,p)=1$ (resp. for every $r\in\N$). 
We then may also adopt the following evocative terminology: a prime $p$ is {\it transparent} to $G$ iff $p\in\corsivoF(G)$, while it is {\it invisible} to $G$ iff $p\in\corsivoF_0(G)$. 
\par 

We can measure the ``transparency'' of a prime $p$ by a kind of $p$-adic valuation defined as follows:
$$
v_{p,G}\defineq \max\{ K\in \N_0 : G(p^{K})=1\},
$$
\par
\noindent
with the convention that : $v_{p,G}\defineq \infty$ \thinspace in case the set on the right-hand side is infinite. Thus a prime $p$ is transparent iff $v_{p,G}\geq 1$ and is invisible iff $v_{p,G}=\infty$. 
Whenever $\corsivoF(G)$ is finite, by abuse of notation : 
$$
P(G)\defineq \prod_{p\in \corsivoF(G)}p.
$$
\par
\noindent
We note (compare Remark 1 in [C3], second version) that 
$$
G:\N \rightarrow \C 
\enspace \hbox{\rm is\enspace a\enspace Ramanujan\enspace coefficient}
\quad \Longrightarrow \quad 
\corsivoF(G) \enspace \hbox{\rm and} \enspace \corsivoF_0(G) \enspace \hbox{\rm are\enspace finite\enspace sets}. 
$$
\par

Since $\corsivoF_0(G)\subseteq \corsivoF(G)$, we distinguish the following three cases, already exposed in [C3] (second version), depending on which of these spectra is empty or not: (recall $G:\N \rightarrow \C$ is multiplicative) 
$$
\eqalign
{
	G \enspace \hbox{\rm is} \enspace \hbox{\stampatello normal} 
	\quad &\definiz \quad
	\corsivoF(G)=\emptyset, 
\cr
	G \enspace \hbox{\rm is} \enspace \hbox{\stampatello sporadic} 
	\quad &\definiz \quad
	 \corsivoF(G)\neq \emptyset \enspace {\rm and } \enspace \corsivoF_0(G)=\emptyset, 
\cr
	G \enspace \hbox{\rm is} \enspace \hbox{\stampatello exotic} 
	\quad &\definiz \quad
	\corsivoF_0(G)\neq \emptyset. 
\cr
}
$$

\par				
\noindent
We can now give our Classification theorem of the multiplicative part of $<\0>$ under some technical hypotheses of convergence. Some convergence hypotheses of this sort cannot be removed completely, as we'll argue in future papers. 
\smallskip
\par
\noindent
{\bf Theorem 1.} {\stampatello Classification of multiplicative Ramanujan coefficients of } $\0$. 
\par
\noindent
{\it Let } $G:\N \rightarrow \C$ {\it be multiplicative. Then one and only one of the following cases happens:} 
\smallskip
\item{A)} $G$ {\it is} {\stampatello normal}. {\it Then, assuming } $\sum_{(r,a)=1}G(r)\mu(r)$ {\it converges pointwise } $\forall a\in \N$, $G$ {\it is a Ramanujan coefficient and} 
$$
G\in <\0>
\enspace \Longleftrightarrow \sum_{q=1}^{\infty}G(q)\mu(q)=0; 
$$
\item{B)} $G$ {\it is} {\stampatello sporadic}. {\it Then, assuming } $\sum_{(r,a)=1}G(r)\mu(r)$ {\it converges pointwise } $\forall a\in \N$, $G$ {\it is a Ramanujan coefficient and} 
$$
G\in <\0>
\enspace \Longleftrightarrow \sum_{(q,P(G))=1}G(q)\mu(q)=0; 
$$
\item{C)} $G$ {\it is} {\stampatello exotic}. {\it Then, assuming } $\sum_{(r,a)=1}G(r)\mu(r)$ {\it converges pointwise } $\forall a\in \N$, $G$ {\it is a Ramanujan coefficient and} 
$$
G\in <\0>. 
$$
\par
The multiplicative Ramanujan coefficient $G_R(q)=1/q$ considered in the Introduction is an example of a normal function, while the function $G_H(q)=1/\varphi(q)$ of Hardy is sporadic, as $\corsivoF(G_H)=\{2\}$ and $\corsivoF_0(G_H)=\emptyset$. 
We will give many more examples of normal and sporadic $G\inzero$ using analytic methods in section 5. 
We shall demonstrate that the examples of Ramanujan and Hardy are not isolated; in a suitable sense, the main feature that makes them belong to the $\0$-cloud is their asymptotic rate of decay to zero along the primes: $G_R(p),G_H(p)\sim 1/p$. 
We refer to Proposition 1 in section 5 for a precise statement. 

An example of exotic $G\inzero$ is given in [C3], and in section 6 we will provide several more examples. 
The most elementary one is the indicator function of the powers of 2, i.e. the multiplicative function $G_{2}(q)$ which is $1$ if $q=2^k$ for some $k\in\N_0$, and is $0$ otherwise (compare Remark 2, $\S6$).  
With a generalization of the concept of exotic functions, we are also able in section 6.2 to produce for the first time some examples (in fact, a very large class) of Ramanujan expansions of $\0$ whose Ramanujan coefficients are not multiplicative. 

Many of the  Ramanujan coefficients in the exotic class, including the example $G_{2}$, give rise to absolutely convergent expansions of $\0$ (sometimes even Ramanujan series with finitely many nonzero terms). 
However, we argue in section 6.1 that there exist exotic Ramanujan coefficients whose associated Ramanujan expansion does not converge absolutely to zero. 
We investigate in more detail the issue of absolute convergence in $\S7$. 

The proof of the theorem above will be given in section 4, after having established some elementary technical lemmas in next section 3.
\par
In a forthcoming paper, we'll discuss the hypotheses of convergence in our classification theorem, in relation with a general problem of summability of multiplicative functions supported over the squarefree numbers, when coprimality conditions are imposed. 
\bigskip
\bigskip
\bigskip
\bigskip
\par
\noindent {\bf 3. Lemmata}

\bigskip

\par
\noindent
We start, at once, with the Lemmas that we'll use in the Proof of our Classification Theorem. 
\medskip
\par
For multiplicative arithmetic functions $G$ we have following Lemma 1, whose proof is immediate. 
\smallskip
\par
\noindent
{\bf Lemma 1.} {\it Let } $G:\N \rightarrow \C$ {\it be } {\stampatello multiplicative, } {\it let } $\corsivoF$  {\it be any finite, non-empty subset of primes and let } $p_1$ {\it be any prime with } $p_1\notin \corsivoF$. {\it Then} 
$$
G(p_1)\neq 1, \enspace \sum_{(r,\{p_1\}\cup \corsivoF)=1}G(r)\mu(r) \enspace {\it converges \enspace and} \enspace 
\sum_{(r,\corsivoF)=1}G(r)\mu(r)=0 
\enspace \Longrightarrow \enspace 
\sum_{(r,\{p_1\}\cup \corsivoF)=1}G(r)\mu(r)=0. 
$$
\par				
\noindent
{\bf Proof.} Simply, for $x\in \N$ large enough in terms of $p_1$ (say, $x>p_1$), abbreviating \enspace $\corsivoF_1 := \{p_1\}\cup \corsivoF$, 
$$
\sum_{{r\le x}\atop {(r,\corsivoF_1)=1}}G(r)\mu(r)=\sum_{{r\le x}\atop {(r,\corsivoF)=1}}G(r)\mu(r)-\sum_{{r\le x}\atop {{(r,\corsivoF)=1}\atop {r\equiv 0(\!\!\bmod p_1)}}}G(r)\mu(r)
=\sum_{{r\le x}\atop {(r,\corsivoF)=1}}G(r)\mu(r)+G(p_1)\sum_{{r\le x/p_1}\atop {(r,\corsivoF_1)=1}}G(r)\mu(r), 
$$
\par
\noindent
whence 
$$
\sum_{{r\le x}\atop {(r,\corsivoF)=1}}G(r)\mu(r)=\sum_{{r\le x}\atop {(r,\corsivoF_1)=1}}G(r)\mu(r)-G(p_1)\sum_{{r\le x/p_1}\atop {(r,\corsivoF_1)=1}}G(r)\mu(r). 
$$
\par
\noindent
The limit as $x\to \infty$ vanishes for LHS (Left Hand Side), while for the right-hand side it is $1-G(p_1)\neq 0$ times the series over $(r,\corsivoF_1)=1$.\hfill $\square$ 

\bigskip

\par
Our second Lemma will be applied for the sporadic and the exotic functions $G$. Proof comes quickly. 
\smallskip
\par
\noindent
{\bf Lemma 2.} {\it Let } $G:\N \rightarrow \C$ {\it be } {\stampatello multiplicative}. {\it Assume there's a prime } $p_0$, {\it with } $G(p_0)=1$. {\it Then} 
$$
p_0\not \vert a 
\enspace {\it and} \enspace 
\sum_{(r,ap_0)=1}G(r)\mu(r)
\enspace {\it converges} \enspace 
\enspace \Longrightarrow \enspace 
\sum_{(r,a)=1}G(r)\mu(r)=0
\enspace \Longrightarrow \enspace 
\sum_{q=1}^{\infty}G(q)c_q(a)=0. 
$$
\par
\noindent
{\bf Proof.} The first implication uses $G(p_0)=1$ (once fixed a large enough $x\in \N$, say $x>p_0$): 
$$
\sum_{{r\le x}\atop {(r,a)=1}}G(r)\mu(r)=\sum_{{r\le x}\atop {{(r,a)=1}\atop {r\equiv 0(\!\!\bmod p_0)}}}G(r)\mu(r)+\sum_{{r\le x}\atop {(r,ap_0)=1}}G(r)\mu(r) 
=-\sum_{{r\le x/p_0}\atop {(r,ap_0)=1}}G(r)\mu(r)+\sum_{{r\le x}\atop {(r,ap_0)=1}}G(r)\mu(r), 
$$
\par
\noindent
whence passing to the limit as $x\to \infty$, 
$$
\sum_{(r,a)=1}G(r)\mu(r)=0, 
$$
\par
\noindent
while the second implication comes from the Proposition in [C3] (second version).\hfill $\square$ 

\vfill
\eject

\par				

\par 
\noindent {\bf 4. Proof of the Classification Theorem}

\bigskip

\par
\noindent
We start with 
\par
\noindent
{\bf Proof of case A).} Follow the proof of Corollary 1 in [C3], second version. \qed 
\medskip
\par
\noindent
Thanks to Lemma 2, we prove case C) immediately. 
\par
\noindent
{\bf Proof of case C).} Apply the Proposition in [C3], second version, getting, $\forall a\in \N$ fixed, 
$$
\sum_{q=1}^{\infty}G(q)c_q(a)=\prod_{p|a}\sum_{K=0}^{v_p(a)}p^K \left(G(p^K)-G(p^{K+1})\right) \,\cdot\, \sum_{(r,a)=1}G(r)\mu(r), 
\leqno{(\ast)_a}
$$
\par
\noindent
for which finite product the condition $(a,\corsivoF_0(G))>1$, implying $\exists p_0\in \P(a)\cap \corsivoF_0(G)$, gives vanishing. 
\par
\noindent
We are in case $(a,\corsivoF_0(G))=1$, otherwise : this time $(\ast)_a$ co-finite factor vanishes, applying Lemma 2. \qed
\medskip
\par
\noindent
Our last, and longest, proof is when $G$ is sporadic. 
\par
\noindent
{\bf Proof of case B).}
\par
\noindent
We quote and prove Corollary 3 of [C3] (second version). For sporadic $G$, notice: \enspace $(r,P(G))=(r,\corsivoF(G))$. The QED (Quod Erat Demonstrandum$=$what was to be shown) will indicate the end of a part of the Proof. 
\smallskip
\par
\noindent
{\bf Corollary 3[C3](2nd ver.).} {\it Let } $G:\N \rightarrow \C$ {\it be } {\stampatello sporadic} {\it and assume that } $\sum_{(r,a)=1}G(r)\mu(r)$ {\it converges pointwise } $\forall a\in \N$. {\it Then } $G$ {\it is a Ramanujan coefficient and} 
$$
\sum_{q=1}^{\infty}G(q)c_q(a)=\0(a) 
\quad 
\Longleftrightarrow 
\quad 
\sum_{(r,\corsivoF(G))=1}G(r)\mu(r)=0. 
$$
\par
\noindent
{\bf Remark 1.} {\it The condition of convergence on the series with $\mu$ immediately implies, from Proposition of} [C3] {\it (second version), that $G$ is a Ramanujan coefficient}.\hfill $\diamond$ 
\medskip
\par
\noindent
{\bf Proof.} We start with \lq \lq $\Rightarrow$\rq \rq, the easiest. 
Since $G$ is a non-exotic Ramanujan coefficient, the $p$-adic valuation $v_{p,G}$ of $G$ is never $\infty$; also, it's non-zero for only finitely many primes: exactly, those in $\corsivoF(G)$. Then the natural number 
$$
a_G\defineq \prod_{p\in \corsivoF(G)} p^{v_{p,G}},
$$
\par
\noindent
is well-defined, and $v_p(a_G)= v_{p,G}$ for every $p\in\P$. 
If we choose \enspace $\corsivoF = \corsivoF(G)$ and \thinspace $a=a_G\in \N$ in the Main Theorem of [C3] (second version), we get the following formula:
$$
0=\sum_{q=1}^{\infty}G(q)c_q(a_G)
=\left(\prod_{p\in \corsivoF(G)}\sum_{K=0}^{v_{p,G}}p^K (G(p^K)-G(p^{K+1}))\right)
 \cdot \sum_{(r,\corsivoF(G))=1}G(r)\mu(r). 
\leqno{(\ast)_G}
$$
\par
\noindent
By the definition of $v_{p,G}$ we have \enspace $G(p^{K})=1$ \enspace for every $K\leq v_{p,G}$ and \enspace $G(p^{v_{p,G}+1})\neq 1$. Then, the {\bf finite factor} (i.e., the finite Euler product) in $(\ast)_G$ becomes 
$$
\prod_{p\in \corsivoF(G)}p^{v_{p,G}}(1-G(p^{v_{p,G}+1}))=a_G \prod_{p\in \corsivoF(G)}(1-G(p^{v_{p,G}+1}))
\neq 0, 
$$
\par
\noindent
forcing the {\bf co-finite factor} (the infinite series complementary to the finite factor) in $(\ast)_G$ to vanish:  
$$
\sum_{(r,\corsivoF(G))=1}G(r)\mu(r)=0.
$$
\par
\noindent
\hfill {\rm QED}($\Rightarrow$) 

\bigskip

\par				
\noindent
The other implication, i,e., \lq \lq $\Leftarrow $\rq \rq, starts from 
$$
\sum_{(r,\corsivoF(G))=1}G(r)\mu(r)=0. 
\leqno{(\ast)}
$$
\par
\noindent
We'll get the vanishing of Ramanujan expansion (LHS in Corollary) little by little, on all $a\in \N$. 
\smallskip
\par
We do it in three steps.
\smallskip
\par
\leftline{STEP 1. }
\par
\noindent
From the Proposition in [C3] (second version) we get, choosing $a\in \N$ with $\P(a)=\corsivoF(G)$, 
$$
\sum_{q=1}^{\infty}G(q)c_q(a)=\left(\prod_{p|a}\sum_{K=0}^{\infty}G(p^K)c_{p^K}(a)\right)\cdot \sum_{(r,\corsivoF(G))=1}G(r)\mu(r)=0,
\quad 
\forall a\in \N \;:\; \P(a)=\corsivoF(G), 
$$
\par
\noindent
from $(\ast)$ above. 
\par
\leftline{STEP 2. }
\par
\noindent
Now we're going to enlarge the set of $a\in \N$ we got in previous step, using Lemma 1.
\par
\noindent
Choosing, in quoted Proposition, $a\in \N$ with $\P(a)=\corsivoF(G)\cup \{p_1\}$, say $p_1\notin \corsivoF(G)$, 
$$
\sum_{q=1}^{\infty}G(q)c_q(a)=\left(\prod_{p|a}\sum_{K=0}^{\infty}G(p^K)c_{p^K}(a)\right)\cdot \sum_{(r,\corsivoF(G)\cup \{p_1\})=1}G(r)\mu(r), 
\quad 
\forall a\in \N \;:\; \P(a)=\corsivoF(G)\cup \{p_1\}, 
$$
\par
\noindent
whence we get, from both $(\ast)$ and Lemma 1, 
$$
\sum_{q=1}^{\infty}G(q)c_q(a)=0, 
\quad 
\forall a\in \N \;:\; \P(a)=\corsivoF(G)\cup \{p_1\}. 
$$
\par
\noindent
Iterating this procedure and joining more different primes, applying inductively Lemma 1, we get 
$$
\sum_{(r,a)=1}G(r)\mu(r)=0, 
\quad 
\forall a\in \N \;:\; a\equiv 0\,(\bmod \; P(G)), 
$$
\par
\noindent
because this series depends only on the squarefree kernel $P(a)$ of $a$ (and not on $a$, actually). This last vanishing can be used again in Proposition of [C3] (second version), so
\par
\leftline{WE PROVED:}
$$
\sum_{q=1}^{\infty}G(q)c_q(a)=0, 
\quad 
\forall a\in \N \;:\; a\equiv 0\,(\bmod \; P(G)). 
$$
\par
\leftline{STEP 3. }
\par
\noindent
We are left with the task to 
\par
\leftline{PROVE:}
$$
\sum_{q=1}^{\infty}G(q)c_q(a)=0, 
\quad 
\forall a\in \N \;:\; a\not \equiv 0(\!\bmod \; P(G)). 
$$
\par
\noindent
This is equivalent to proving : 
$$
\exists p_0\in \corsivoF(G), p_0 \not \vert a
\enspace \Longrightarrow \enspace 
\sum_{q=1}^{\infty}G(q)c_q(a)=0. 
$$
\par
\noindent
Lemma 2 is just what we need to conclude.\hfill $\square$ 

\vfill
\eject

\par				
\noindent {\bf 5. Examples of Ramanujan expansions of $\0$ via analytic methods}
\bigskip
\par
\noindent
As we remarked in the introduction, the Ramanujan expansions of $\0$ given by Ramanujan and Hardy, corresponding respectively to the multiplicative functions $G_R(q)=1/q$ and $G_H(q)=1/\varphi(q)$, are related in an essential way to the Prime Number Theorem. 
A close inspection of Hardy's proof [H, section 8] reveals that this fact does not depend much on the special choice of the multiplicative function $G=G_R$ or $G=G_H$, but only on the asymptotic behaviour of $G$. With this intuition, we are able to exhibit many more examples of Ramanujan expansions of $\0$. 

\bigskip
\par
\noindent
{\bf Proposition 1.} {\it Let $G$ be multiplicative such that  $G(q)=O( 1 / q )$ for $q$ squarefree and 
$$
G(p)= {1 \over p}+O(p^{-1-\alpha})
$$ 
\par
\noindent
for some $\alpha >0$ as $p\to\infty$ along the primes. Then $G$ is a Ramanujan coefficient and $G\inzero$.} 
\bigskip

The proof of this Proposition relies on the following Lemma, which we prove by analytic means.

\bigskip
\par
\noindent
{\bf Lemma 3.} {\it Let $\alpha$ and $G$ be as in the Proposition above. Then} 
$$
\sum_{q=1}^{\infty} G(q)\mu(q) = 0.
$$ 
\bigskip
\par
\noindent
{\it Proof of Lemma 3.} First we define the Dirichlet series $$
f(s) := \sum_{q=1}^{\infty} {G(q)\mu(q) \over q^{s-1}},
$$ 
\par
\noindent
which is absolutely convergent for $\sigma \defineq Re(s)>1$ and is equal to the Euler product
$$
f(s)=\prod_{p\in\P} ( 1 - G(p) p^{-s+1} )
$$ 
on the same half plane. We recall that the Riemann zeta function is given by the Euler product
$$
\zeta(s) = \prod_{p\in \P} (1- p^{-s})^{-1}
$$
\par
\noindent
for $\sigma >1$. Then, if we write $G(p)=(1+\epsilon(p))p^{-1}$ for some $\epsilon(p)=O(p^{-\alpha})$ we get that, say, 
$$
h(s):=f(s)\zeta(s) = \prod_{p\in \P} { 1-(1+\epsilon(p)) p^{-s} \over 1-p^{-s} }.
$$ 
Since every factor in this product is of the form 
$$
1-\epsilon(p)p^{-s} + O(p^{-2\sigma -\alpha}) = 1+O(p^{-\sigma -\alpha}), 
$$
\par
\noindent
we have on the half plane $\sigma>1$ that 
$$
f(s) = h(s) \cdot {1\over \zeta(s) }, 
$$

\vfill
\eject

\par				
\noindent
for some Dirichlet series $h(s)$ that extends analytically to the half plane $\sigma > 1 -\alpha$. 
Moreover, by the analytic version of the Prime Number Theorem, we know that the inverse of the zeta function extends by continuity to the closed half plane $\sigma \geq 1$, and so does $f(s)$. Then, by the Theorem of [N], we get that $\sum_{q=1}^{\infty} G(q)\mu(q)$ converges and is equal to the value of $\lim_{s \to 1} f(s)$. Since the function $h(s)$ is regular at $s=1$ and the zeta function $\zeta(s)$ has a pole there, we get that this limit is equal to zero. \qed

\bigskip
\par
\noindent
{\it Proof of Proposition 1.} Fix any $a\in\N$. We may apply the previous Lemma to the multiplicative function $G \,\cdot\, \1_{(q,a)=1}$, which is equal to $G$ on the numbers coprime to $a$ and is zero elsewhere. 
Then we get
$$
\sum_{(q,a)=1} G(q)\mu(q)=0.
$$ 
\par
\noindent
In particular, this series converges for every $a$. Since $G$ is multiplicative, the Proposition of [C3] (second version) applies and we get $G\inzero$. \qed

\bigskip
\bigskip
\bigskip
\par
\noindent {\bf 6. Examples of exotic and weakly-exotic Ramanujan coefficients}
\bigskip
\par
\noindent 
It was already observed by the first author in [C3] (second version) that there exist exotic Ramanujan coefficients \thinspace $G\in\,<\0>$ \thinspace for which the corresponding Ramanujan series \thinspace $\sum_{q=1}^\infty G(q)c_q(a)$ \thinspace converges absolutely to zero, for every \hfil $a\in\N$. In fact, there are lots of absolutely convergent Ramanujan expansions of this kind, as the following construction shows. 
\bigskip
\par
\noindent
{\bf Proposition 2.} {\it Let $p_0$ be a prime number and let $G$ be a multiplicative function such that $G(p_0^k)=1$ for all $k\in\N_0$ and such that $\sum_{p\in \P}^{} |G(p)|<\infty.$ Then $G$ is an exotic Ramanujan coefficient and} 
$$
\sum_{q=1}^{\infty} G(q)c_q(a) 
$$
\par
\noindent
{\it converges absolutely to $0$ for all $a\in\N$.} 
\medskip

The proof of this result is postponed to the next section (see Lemma 5 in $\S7.1$), where we will discuss in more detail the absolutely convergent Ramanujan expansions with multiplicative coefficients. 
In fact, we will also be able to show that all multiplicative coefficients $G\inzero$ with absolutely convergent Ramanujan expansions come from Proposition 2, see Theorem 3. 
\medskip

\par
\noindent
{\bf Remark 2.} {\it A trivial application of Proposition 2, actually, gives the simplest examples of} {\stampatello exotic} {\it function, namely, once fixed any prime } $\tilde{p}$, {\it the function defined as}
$$
G_{\tilde{p}}(q)=1,
\quad
\hbox{\it if} \enspace q=\tilde{p}^K,
\quad \forall K\ge 0,
$$
{\it and } $0$ {\it otherwise, is clearly multiplicative and}
$$
\sum_{q=1}^{\infty}G_{\tilde{p}}(q)c_q(a)=\sum_{K=0}^{v_{\tilde{p}}(a)}\tilde{p}^K\,(G_{\tilde{p}}(\tilde{p}^K)-G_{\tilde{p}}(\tilde{p}^{K+1}))
=0,
\quad \forall a\in \N,
$$
\par
\noindent
{\it see the Main Lemma in second version of } [C3], {\it with absolute convergence}: 
$$
\sum_{q=1}^{\infty}\left|G_{\tilde{p}}(q)c_q{a}\right|=\sum_{K=0}^{\infty}\left|c_{\tilde{p}^K}(a)\right|
=\sum_{K=0}^{v_{\tilde{p}}(a)+1}\left|c_{\tilde{p}^K}(a)\right|
<\infty,
\quad \forall a\in \N,
$$
\par
\noindent
{\it again by quoted Lemma; the condition in Proposition 2 is trivial}: 
$$
\sum_{p\in \P}\left|G_{\tilde{p}}(p)\right|=\left|G_{\tilde{p}}(\tilde{p})\right|
=1.
$$
\par
\hfill $\diamond$

\vfill
\eject

\par				
\noindent{\bf 6.1 A non-absolutely convergent Ramanujan expansion with exotic coefficients}
\bigskip
\par
\noindent
Despite this large class of examples, we would like to remark that not all exotic Ramanujan coefficients give rise to absolutely convergent Ramanujan series. Here is an example. (It's a small Lemma : we call it \lq \lq Fact\rq \rq) 

\bigskip
\par
\noindent
{\bf Fact (Example).} {\it Let $p_0$ be a prime number and let  $G_0$ be a multiplicative function such that $G_0(p_0^k)=1$ for all $k\in\N$ and $G_0(p)=p^{-1}$ for all primes $p\neq p_0$. Then for all $a\in\N$ the Ramanujan expansion} 
$$
\corsivoR_0(a):=\sum_{q=1}^\infty G_0(q) c_q(a) 
$$
\par
\noindent
{\it converges pointwise to zero, but it does not converge absolutely}.\hfill $\diamond$
\bigskip
\par
\noindent
{\it Proof.} The gist of the argument is that for every $b\in\N$ the following series converges to $0$
$$
\sum_{(q,b)=1} {1\over q} \,\mu(q)=0. 
\leqno{(PNT)_b}
$$
\par
\noindent
This is a consequence of the Prime Number Theorem (PNT) in arithmetic progressions. The rate of convergence to zero of series of this type was examined for instance by Ramar\'e [Ra]. 
\par
\noindent
Taking $G=G_0$, $a\not\equiv 0(\bmod \  p_0)$ and \enspace $\sum_{(r,ap_0)=1}G_0(r)\mu(r)=0$ \enspace from $(PNT)_{ap_0}$, Lemma 2 gives $\corsivoR_0(a)=0$, $\forall a\not\equiv 0(\!\bmod \; p_0)$; while \hfil $\forall a\equiv 0(\!\bmod \; p_0)$, from Proposition [C3], second version, \hfil $(PNT)_a$ $\Rightarrow$ $\corsivoR_0(a)=0$, too. 
\par
\noindent
The fact that the convergence is not absolute follows from the divergence of $\sum_{p\in \P}|G_0(p)|$ and by the criterion of absolute convergence in section 7.1. \qed

\bigskip
\par
\noindent{\bf 6.2 Another big class of Ramanujan coefficients in $<0>$}

\bigskip

\par
\noindent
We define $G:\N \rightarrow \C$ to be {\stampatello weakly exotic} iff 
$$
\exists p_0\in \P \enspace : \enspace G(p_0^K r)=G(r), \enspace \forall K\in \N_0, \enspace \forall r\in \N, (r,p_0)=1 
\leqno{(\ast)_{we}}
$$
\par
\noindent
and we call any prime satisfying $(\ast)_{we}$ above an \lq \lq {\it invisible}\rq \rq \enspace prime for $G$; thus, generalizing the same concept from exotic case, but see that now we do not require $G:\N \rightarrow \C$ to be multiplicative. However, notice that  :  $G$ exotic $\Rightarrow $ $G$ weakly exotic (of course!). Notwithstanding this, we need, for the next result to hold, a stronger convergence condition, w.r.t. the one required by exotic functions, compare case $C)$ of our Classification (Theorem 1 in section 2). 
\par
\noindent
{\bf Theorem 2.} {\stampatello Weakly-exotic functions satisfying a convergence condition are in $\0-$cloud}. 
\par
\noindent
{\it Let } $G:\N \rightarrow \C$ {\it be } {\stampatello weakly-exotic} {\it and, for at least one invisible prime } $p_0$, {\it assume that} 
$$
\sum_{(r,p_0)=1}G(r)c_r(a)
\quad {\it converges}\enspace {\it pointwise}\quad
\forall a\in \N. 
$$ 
\par 
\noindent
{\it Then } $G$ {\it is a Ramanujan coefficient and } \enspace 
$$
\sum_{q=1}^{\infty}G(q)c_q(a)=\0(a). 
$$ 
\par
\noindent
{\bf Proof.} Let's fix $a\in \N$ and choose $Q\in \N$, large enough in terms of $a$ and $p_0$ above (say, $Q>p_0^{v_{p_0}(a)+1}$) : 
$$
\sum_{q\le Q}G(q)c_q(a)=\sum_{K=0}^{\infty}\sum_{{(r,p_0)=1}\atop {r\le Q/p_0^K}}G(p_0^K r)c_{p_0^K}(a)c_r(a)
=\sum_{K=0}^{v_{p_0}(a)+1}c_{p_0^K}(a)\sum_{{(r,p_0)=1}\atop {r\le Q/p_0^K}}G(r)c_r(a), 
$$
\par				
\noindent
because : $c_{p_0^K}(a)=0$, $\forall K>v_{p_0}(a)+1$, see Fact 1 before Main Lemma in second version of [C3]; whence, passing to the limit over $Q\to \infty$, we get $G$ is a Ramanujan coefficient and, since : 
$$
\sum_{K=0}^{v_{p_0}(a)+1}c_{p_0^K}(a)=\sum_{K=0}^{v_{p_0}(a)}\varphi(p_0^K)-p^{v_{p_0}(a)}
=0, 
$$
\par
\noindent
compare quoted Lemma, we also get $G\in<\0>$.\hfill $\square$ 
\medskip
\par
\noindent
Thus the hypothesis \lq \lq $G$ multiplicative\rq \rq, actually, is not necessary for the weakly exotic $G$. 
\medskip
\par
\noindent
The Class of weakly exotic $G:\N \rightarrow \C$ seems to be the widest we know, at the moment, in the cloud of $\0$. 

\bigskip
\bigskip

\par
\noindent {\bf 7. Absolutely convergent Ramanujan expansions}

\bigskip

\par
\noindent
Since lot of research on Ramanujan expansions has been carried out in the hypothesis of absolute convergence (see for instance the theorems of Wintner [W] and Delange [De] or the discussions in the surveys [C1], [L2]), it is natural to study those Ramanujan expansions of $\0$ which are absolutely convergent. 
The classical examples of Hardy and Ramanujan are known to converge pointwise but not absolutely, while Proposition 2 in the previous section shows that there are many exotic coefficients which imply absolute convergence. 
In this section we examine this topic in detail, thus completing a discussion in [C3].
We show that all Ramanujan expansions of $\0$ with {\stampatello normal} and {\stampatello sporadic} multiplicative Ramanujan coefficients are necessarily {\it not} absolutely convergent. 
We also show that all absolutely convergent Ramanujan expansions with exotic coefficients are constructed as in Proposition 2 (section 6). 


\bigskip
\bigskip
\bigskip

\par
\noindent{\bf 7.1 A criterion for absolute convergence}
\bigskip
\par
\noindent
The problem of absolute convergence for Ramanujan expansions with multiplicative coefficients is easily solved using the theory of finite-cofinite Euler product decomposition developed in [C3]. Indeed, we have the following simple criterion.

\bigskip
\par
\noindent
{\bf Lemma 4.} {\it Let $G$ be a multiplicative function. Then the following are equivalent:}
\smallskip
i) \enspace $\sum_{q=1}^{\infty} G(q)c_q(a)$ {\it is absolutely convergent for every $a\in\N$;}
\smallskip
ii)\enspace $\sum_{q=1}^{\infty} G(q)c_q(a)$ {\it is absolutely convergent for some $a\in\N$;}
\smallskip
iii) \enspace$\sum_{q=1}^{\infty} G(q)\mu(q)$ {\it is absolutely convergent;}
\smallskip
iv) \enspace $\sum_{p\in\P} |G(p)|<\infty$.
\bigskip
\par
\noindent
{\it Proof.} The implications (i) $\Rightarrow$ (ii) and (iii) $\Rightarrow$ (iv) are obvious; using $|c_q(a)|\ge \mu^2(q)$, $\forall a,q\in \N$, we get (ii) $\Rightarrow$ (iii). 
So, now let us assume (iv).
The assertion (iii) is equivalent to the summability of $|G|$ over the squarefree numbers. 
Hence, (iii) follows from (iv), using the absolute summability over the prime numbers, thanks to the multiplicativity of $G$, $G(p)\to 0$ and the inequality $\log(1+|z|)\leq |z|$ ($\forall z\in \C$ with $|z|<1$) : 
$$
\sum_{q\  {\rm squarefree}} |G(q)| = \prod_{p\in \P} (1+|G(p)|) \ll  \exp\left(\sum_{p\in\P} |G(p)|\right)<\infty.
$$ 
\par
\noindent
Now, let us prove (i) and (ii) assuming (iii). 
For each $a\in\N$ we recall that $c_q(a)=\mu(q)$ if $(q,a)=1$ and that $c_q(a)=0$ whenever $v_p(q)>v_p(a)+1$ for some prime number $p$. Together with the multiplicativity of $G$ and of the Ramanujan coefficients, we get the following formula, which is analogous to the finite-cofinite Euler product formula in the Main Lemma of [C3]: 
$$
\sum_{q=1}^{\infty} |G(q)c_q(a)| = \sum_{d\mid a P(a)} |G(d)c_d(a)| \sum_{(r,a)=1} |G(r)\mu(r)|
$$
\par				
\noindent
where we recall that $P(a)$ is the squarefree kernel of  $a$, so that $a P(a) = \prod_{p\mid a} p^{v_{p}(a)+1}$.  
Since the first factor is a finite sum and the second factor is a sum over squarefree numbers, we deduce that the Ramanujan series $\sum_{q=1}^{\infty} G(q)c_q(a)$ is absolutely summable. Since $a\in\N$ is arbitrary, the result follows. \qed

\bigskip

A direct consequence of this criterion of absolute convergence is that the Proposition 2 (section 6) is true and, moreover, it captures the structure of all absolutely convergent Ramanujan expansions with exotic coefficients. 

\bigskip
\par
\noindent{\bf 7.2 Only pointwise convergence for {\stampatello normal} and {\stampatello sporadic} coefficients}
\bigskip
\par
\noindent
We now examine the problem of absolute convergence of non-exotic Ramanujan coefficients in the cloud of $\0$.
We prove that absolute convergence is never achieved in this case. 
\bigskip
\par
\noindent
{\bf Proposition 3.} {\it Let $G\in\,<\0>$ be a normal or sporadic Ramanujan coefficient. Then, for each $a\in\N$, the Ramanujan series} 
$$
\sum_{q=1}^\infty G(q) c_q(a),
$$
\par
\noindent
{\it does not converge absolutely (but it converges pointwise to $0$).} 
\bigskip
\par 
\noindent
{\it Proof.} Let $G\in\,<\0>$ be a normal or sporadic Ramanujan coefficient such that the series
$$
\sum_{q=1}^\infty G(q) c_q(a)
$$
\par
\noindent
converges absolutely for some $a\in\N$. We are going to derive a contradiction from this assumption. 
By the criterion of absolute convergence in the previous section (Lemma $4$), we have that
$$
\sum_{r \  {\rm squarefree}} |G(r)|<\infty
\quad \quad {\rm and} \quad \quad 
\sum_{p\in \P} |G(p)| <\infty.
$$
We deduce for each $b\in\N$ that the series
$$
\sum_{(r,b)=1} G(r) \mu(r)
$$
converges and (since $G$ is multiplicative) is equal to the (also convergent) Euler product
$$
\prod_{{p\ {\rm prime}} \atop {(p,b)=1}} (1-G(p)).
$$ 
Since $G\in\,<\0>$, parts $A)$ and $B)$ of Theorem 1 imply that
$$
\sum_{(r,P(G))=1} G(r)\mu(r) = 0,
$$
where  $P(G)=1$ in case $G$ is normal. 
However we also have (compare Property 2 [C3], second version)
$$
\prod_{{p\ {\rm prime}} \atop {(p,P(G))=1}} (1-G(p)) \neq 0,
$$  
because $\sum_{p\in \P} |G(p)|<\infty$  and because $G(p)\neq 1$ for all primes appearing in this Euler product, by definition of $P(G)$. 
This is the  required contradiction. \qed

\bigskip
We summarize our findings as follows.
\bigskip
\par
\noindent
{\bf Theorem 3.} {\stampatello Characterization of absolute convergence in multiplicative part of $\0-$cloud}. 
\par
\noindent
{\it Let $G$ be a multiplicative Ramanujan coefficient in the $\0$-cloud. We have that $\sum_{q=1}^{\infty} G(q) c_q(a)$ converges absolutely for some (hence all) $a\in\N$ if, and only if:} 
$$
\sum_{p\in \P} |G(p)| <\infty      
\quad\quad 
\hbox{\stampatello and}
\quad\quad
\corsivoF_0(G)\neq \emptyset. 
$$

\vfill
\eject

\par				
\centerline{\stampatello Bibliography}

\bigskip

\item{[Ca]}  R. Carmichael, {\sl Expansions of arithmetical functions in infinite series}, Proc. London Math.
Society, {\bf 34} (1932), 1-26.
\smallskip
\item{[C1]} G. Coppola, {\sl A map of Ramanujan expansions}, ArXiV:1712.02970v2. (Second Version) 
\smallskip
\item{[C2]} G. Coppola,{\sl A smooth shift approach for a Ramanujan expansion}, ArXiV:1901.01584v3.(Third Version) 
\smallskip
\item{[C3]} G. Coppola, {\sl Finite and infinite Euler products of Ramanujan expansions}, ArXiV:1910.14640v2 (Second Version) 
\smallskip
\item{[C4]} G. Coppola, {\sl Recent results on Ramanujan expansions with applications to correlations}, to appear on Rend. Semin. Mat. Univ. Politec. Torino 
\smallskip
\item{[CMS]} G. Coppola, M. Ram Murty and B. Saha, {\sl Finite Ramanujan expansions and shifted convolution sums of arithmetical functions},  J. Number Theory {\bf 174} (2017), 78--92. 
\smallskip
\item{[CM]} G. Coppola and M. Ram Murty, {\sl Finite Ramanujan expansions and shifted convolution sums of arithmetical functions, II}, J. Number Theory {\bf 185} (2018), 16--47. 
\smallskip
\item{[D]} H. Davenport, {\sl Multiplicative Number Theory}, 3rd ed., GTM {\bf 74}, Springer, New York, 2000. 
\smallskip
\item{[De]} H. Delange, {\sl On Ramanujan expansions of certain arithmetical functions}, Acta Arith. {\bf 31}(1976), 259--270.
\smallskip
\item{[H]} G.H. Hardy, {\sl Note on Ramanujan's trigonometrical function $c_q(n)$ and certain series of arithmetical functions}, Proc. Cambridge Phil. Soc., {\bf 20} (1921), 263--271. 
\smallskip
\item{[H\"o]} O. H\"older, {\sl Zur Theorie der Kreisteilungsgleichung $K_m (x)=0$}, Prace Mat. Fiz. {\bf 43} (1936), 13--23.
\smallskip
\item{[Hoo]} C. Hooley, {\sl A note on square-free numbers in arithmetic progressions}, Bull. London Math. Soc. {\bf 7} (1975), 133-138.
\smallskip
\item{[K]} J.C. Kluyver. {\sl Some formulae concerning the integers less than $n$ and prime to $n$}, Proceedings of the Royal Netherlands Academy of Arts and Sciences (KNAW), {\bf 9} (1) (1906), 408--414. 
\smallskip
\item{[L1]} L.G. Lucht, {\sl Ramanujan expansions revisited}, Arch. Math. {\bf 64} (1995), 121--128. 
\smallskip
\item{[L2]} L.G. Lucht, {\sl A survey of Ramanujan expansions}, International Journal of Number Theory, {\bf 6} (2010), 1785-1799.
\smallskip
\item{[M]} M. Ram Murty, {\sl Ramanujan series for arithmetical functions}, Hardy-Ramanujan J. {\bf 36} (2013), 21--33. Available online 
\smallskip
\item{[N]} D. Newman, {\sl Simple analytic proof of the prime number theorem}, Amer. Math. Monthly {\bf 87} (1980), no. 9, 693--696.
\smallskip
\item{[R]} S. Ramanujan, {\sl On certain trigonometrical sums and their application to the theory of numbers}, Transactions Cambr. Phil. Soc. {\bf 22}  (1918), 259--276.
\smallskip
\item{[Ra]}  O. Ramar\'e, Explicit estimates on the summatory functions of the M\"obius function with coprimality restrictions, Acta Arith. {\bf 165}  (2014), 1--10.
\smallskip
\item{[ScSp]} W. Schwarz and J. Spilker, {\sl Arithmetical Functions}, Cambridge University Press, 1994.
\smallskip
\item{[T]} G. Tenenbaum, {\sl Introduction to Analytic and Probabilistic Number Theory}, Cambridge Studies in Advanced Mathematics, {\bf 46}, Cambridge University Press, 1995. 
\smallskip
\item{[W]} A. Wintner, {\sl Eratosthenian averages}, Waverly Press, Baltimore, MD, 1943. 

\bigskip
\bigskip
\bigskip

\leftline{\tt Giovanni Coppola - Universit\`{a} degli Studi di Salerno (affiliation)}
\leftline{\tt Home address : Via Partenio 12 - 83100, Avellino (AV) - ITALY}
\leftline{\tt e-mail : giovanni.coppola@unina.it}
\leftline{\tt e-page : www.giovannicoppola.name}
\leftline{\tt e-site : www.researchgate.net}

\bigskip
\bigskip

\leftline{\tt Luca Ghidelli - Universit\`{a} degli Studi di Genova (affiliation)}
\leftline{\tt Address : Via Dodecaneso 35 - CAP 16146, Genova (GE) - ITALY}
\leftline{\tt e-mail : luca.ghidelli@uottawa.ca}
\leftline{\tt e-page : math-lucaghidelli.site}

\bye